\input amstex
\documentstyle{amsppt}
\magnification\magstep1         
\baselineskip=10pt
\parindent=0.2truein
\pagewidth{6.6truein}
\pageheight{9.0truein}
\topmatter
\endsubjclass

\title  Remarks on Ruelle Operator and\\
invariant line fields problem II
\endtitle
\author  Peter M. Makienko
\endauthor
\abstract{Let $ R $ be rational map. Critical point $ c $ is called {\it summable} if  series $ \sum_i\frac{1}{(R^i)'(R(c))} $ is absolutely convergent.  Under some topological condition on postcritical set we prove that $ R $ can not be structurally stable if summable critical point $ c \in J(R). $}
\endabstract
\thanks
This work was partially supported by proyecto {\bf CONACyT}
{\it Sistemas Din\'amicos} 27958E from Mexico
\endthanks

\address{Permanent addresses:
P. Makienko, Instituto de Matematicas,\newline 
Av. de Universidad s/N., Col. Lomas de Chamilpa, C.P. 62210, \newline Cuernavaca, Morelos, Mexico and \newline Institute for Applied Mathematics, \newline  9 Shevchenko str.,
Khabarovsk, Russia}
\endaddress
\thanks{This work has been partially supported by the Russian Fund of Basic
Researches, Grant 99-01-01006.}
\endthanks

\endtopmatter
\document
\heading{\bf Introduction and main statement}\endheading
\par This work is continuation and generalization of works \cite{\bf Lev}, \cite{\bf Mak} and \cite{\bf Mak1}. In summable case, see definition below,  the idea of G. Levin on construction of fixed point for Ruelle operator allows to remove certain cumbersome conditions of \cite{\bf Mak}. 
\par The main aim of this work is to prove theorem A below.
Let $ c_i $ be critical points of $ R $ and  $ Pc(R) = \overline{\cup_i\cup_{n \geq 0}R^n(R(c_i))}$ be the postcritical set.
\proclaim{Definition} A  point $ a $ of for a given rational map $ R $ is called "{\rm summable}" if and only if either \hfill\newline
{\rm 1)}the set $ X_a = \overline{\{\cup_nR^n(R(a))\}} $ is bounded and the series
$$\sum_{i = 0}\frac{1}{(R^i)\prime(R(a))}  
$$ 
is absolutely convergent or\hfill\newline
{\rm 2)}the set $ X_a $ is unbounded and the series
$$
\sum_{i = 0}\frac{1}{(R^i)\prime(R(a))}  \text{ and } \sum_{i = 0}\frac{\vert R^n(R(a)\vert\vert \ln\vert R^n(R(a)\vert\vert}{(R^i)\prime(R(a))}
$$ 
are absolutely convergent.
\endproclaim

\par Note that the property of "{\it summability}" in definition above is not invariant respect to conjugation by Mobius maps. Indeed let  the set $ X_c $ be unbounded and let $ h $ be a Mobius map mapping $ X_c $ into complex plane. Then summability point $ h(c) $ for map $ h\circ R\circ h^{-1} $ is equivalent to absolutely convergence of the series
$ \sum_{i = 0}\frac{(R^n(R(c)))^2}{(R^i)\prime(R(c))}  $
for initial map $ R. $ Hence condition 2) in definition above a little bit weaker then "conjugated version " of condition 1).

\proclaim{Definition} Let $ X $ be the space of rational maps $ R $ fixing $ 0, 1, \infty $(if $ F(R) \neq \emptyset, $ then we assume $\infty \in F(R)$)  with summable critical point $ c \in J(R) $ and either
\roster 
\item $ c \notin X_c $ or
\item $ X_c $ does not separate plane or
\item $ m(X_c) = 0, $ where $ m $ is the Lebesgue measure or
\item $ c \in \partial D \subset J(R),$ where
$ D $ is some component of $ F(R).$ 
\endroster
\par Note that last case includes the maps with completely invariant domain.
\endproclaim

\proclaim{\bf Definition} The next set is called space of quasiconformal
deformations of a given rational map $ R $ and denoted by $ qc(R). $
$$
\aligned qc(R) = \bigl\{
& F \in {\Bbb C}P^{2d+1}: \text{ there is a quasiconformal
automorphism $ h_F $ of the Riemann}\\
&\text{ sphere } \overline{\Bbb C} \text{ such that }
 F = h_F \circ R \circ {h^{-1}_F}\bigr\}{\big/}PSL_2(\Bbb C).
\endaligned 
$$
If $ dim(qc(R)) = 2\deg(R) - 2, $ then $ R $ is structurally stable map
\endproclaim
Now we can formulate the main result of this work.
\proclaim{Theorem A} Let $ R \in X. $
Then $ R $ is not structurally stable map or is {\rm unstable map.}
\endproclaim
\par A. Avila (personal communication and see also {\cite{\bf Lev} and \cite{\bf Mak} for polynomials degree two) can prove Theorem A for polynomials under addition conditions. That is "{\it If for polynomial $ P $ the critical point $ c \in J(P) $ is summable and  $ \sum_{i = 0}\frac{1}{(R^i)\prime(R(c))} \neq 0, $ then $ P $ is unstable.} His ideas using polynomial-like stuff and non every map with completely invariant domain generates a polynomial-like map. 
\proclaim{Corollary A} Let $ R $ be a rational map with summable critical point $ c \in J(R).$ If the sum $ \sum_{i = 0}\frac{1}{(R^i)\prime(R(c))} \neq 0, $ then $ R $ is unstable map.
\endproclaim

	\par The next result (theorem B) is an application of ideas and arguments of Theorem A. 
\proclaim{Definition} Denote by  $ W $ the following subset of rational maps.
\roster 
\item There is no parabolic points for $ R \in W $ and
\item All critical point are simple (that is $ R''(c) \neq 0$) and the forward orbit of any critical point $ c $ is infinite and does not intersect the forward orbit of any other critical point.
\endroster 
\endproclaim
\proclaim{Definition} We call  a rational map $ R $ {\rm summable} if all critical points  belonging to Julia set are summable. 
\endproclaim 
\proclaim{Theorem B} Let $ R \in W $ be summable rational map with completely invariant domain. Then there exists no invariant line fields on $ J(R).$
\endproclaim
\par This result is not new (see results of H. Bruin and S. van Strien {\cite{\bf BS}} and J. Rivera-Letelier {\cite{\bf RL}) and is not so strong like \cite{\bf BS} and \cite{\bf RL}. We just give another approach to invariant line fields problem.
\par  In  forthcoming papers we will discuss application of proposed approach to exponential map and entire functions without asymptotic values.
\subheading{Acknowledgement} We thank to Gena Levin for useful discussion, Alex Eremenko and seminar "Dynamical Systems" IMUNAM at Cuernavaca, Mexico  for useful discussion this work.
\heading{\bf Proof of Theorem A}\endheading
	\par Our aim is to show that under of theorem $ dim(qc(R)) < 2deg(R) - 2, $ then we obtain contradiction with $ dim(qc_J(R)) = 2 deg(R) - 2. $

\par Assume inverse and start with a structurally stable map $ R \in X. $ We always can assume that
\roster
\item  $ R''(c) \neq 0 $ for any critical point $ c. $
\item All periodic components are attractive and the corresponding Riemann surface $ S(R)$ is the union  of punctured torii. Then
$ dim(T(R)) = dim(qc(R)) = 2deg(R)- 2,$
where $ T(R) $ is Teichmuller space for $ R $ (for details see \cite{\bf MS} and \cite{\bf S}.) 
\endroster
\par Now we accumulate some definitions and notations.
\proclaim{Definition} Consider the following space
$$\multline
K(R) = \biggl\{
g\in Rat_{deg(R)}; \text{ there exists a qc-map } f:\overline{\Bbb C}\rightarrow\overline{\Bbb C} \text{ fixing points }\\ 0, 1, \infty
\text{ and } g =f\circ R\circ f^{-1}\biggr\},
\endmultline
$$
then by definition $ K(R)_{/Psl(2, {\Bbb C})}\cong qc(R) $ and  $ dim(K(R)) = dim(qc(R)).$
\endproclaim
\proclaim{Definition Space $H^1(R)$} Let $ h $ be a germ of holomorphic functions at $ 0 \in {\Bbb C} $ with values in $ K(R) $ such that $ h(0) = R. $ Two germs $ h_1 $ and $ h_2 $ we call equivalent ($ h_1 \sim h_2$) iff 
$$
\frac{\partial}{\partial\lambda}h_1(\lambda)_{\vert\lambda = 0} = \frac{\partial}{\partial\lambda}h_2(\lambda)_{\vert\lambda = 0}.
$$
Then we define
$$ H^1(R) = {\biggl\{\text{ Germs like above} \biggr\}}_{{\bigg /}\sim}
$$
or we can write $ h(\lambda) = R(z) + \lambda G_h(z) + ... $, then

$$ 
H^1(R) = \{G_h(z)\}
$$ 
Or if $ R(z) = z\frac{P(z)}{Q(z)} $ with $ deg(P) = deg(R) - 1, deg(Q) = deg(R) - 1 $ and $ P(1) = Q(1),$ and $ h(\lambda) = R_\lambda(z) = z\frac{P + \lambda P_1 + ...}{Q + \lambda Q_1 + ...} = R(z) + \lambda z\frac{P_1Q - Q_1P}{Q^2} + ... $ then 
$$ 
\multline H^1(R) = \biggl\{z\frac{P_1Q - Q_1P}{Q^2}, \text{ where } Q_1(1) = P_1(1), deg(Q_1) \leq deg(R) - 2, \\ 
 deg(P_1) \leq \deg(R) - 1\biggr\},
\endmultline
$$

and $ dim(H^1(R)) = 2deg(R) - 2.$
\endproclaim
\subheading{Bers map}
Let $ \phi \in L_\infty(\Bbb C) $ and $ B_R(\phi) = \phi(R)\frac{{\overline R'}}{R'} : L_\infty(\Bbb C) \rightarrow L_\infty(\Bbb C) $ be Beltrami operator. Then the open unit ball $ B $ of the space $ Fix(R) \subset L_\infty(\Bbb C) $ of fixed points of $ B_R $ called {\it the space of Beltrami differentials for $ R $} and describe all $qc-$ deformations of $ R. $
\par Let $ \mu \in Fix(R), $ then for any $ \lambda $ with $ \vert\lambda\vert < \frac{1}{\Vert\mu\Vert} $ the element $ \mu_\lambda = \lambda\mu \in B. $ Let $ f_\lambda $ be qc-maps corresponding to Beltrami differentials $ \mu_\lambda $ with $ f_\lambda(0, 1, \infty) = (0, 1, \infty).$ Then the map
$$ \lambda \rightarrow R_\lambda = f_\lambda\circ R\circ f_\lambda^{-1} \in K(R) 
$$
is a conformal map. Let $ R_\lambda(z) = R(z) + \lambda G_\mu(z) + ..., $
then differentiation respect to $ \lambda $ in point $ \lambda = 0$ gives the following equation
$$ F_\mu(R(z)) - R'(z)F_\mu(z) = G_\mu(z),
$$
where $ F_\mu(z) = \frac{\partial}{\partial\lambda}f_\lambda(z)_{\vert\lambda = 0}.$
\proclaim{Remark 1} Due to $qc$ map theory (see for example \cite{\bf Krush}) for any $ \mu \in L_\infty(\Bbb C) $ with $\Vert\mu\Vert_\infty < \epsilon $ and small $ \epsilon$ there exists the following formula for qc-map $ f_\mu $ fixing $ 0, 1, \infty.$
$$ 
f_\mu(z) = z -\frac{z(z - 1)}{\pi}\iint_{\Bbb C}\frac{\mu}{\xi(\xi - 1)(\xi - z)} + C(\epsilon, R)\Vert\mu\Vert_\infty^2,
$$
where $ \vert z\vert < R $ and $ C(\epsilon, R) $ is constant does not depending on $ \mu.$ Then
$$
F_\mu(z) = \frac{\partial f_\lambda}{\partial\lambda}_{\vert\lambda = 0} = - \frac{z(z - 1)}{\pi}\iint_{\Bbb C}\frac{\mu}{\xi(\xi - 1)(\xi - z)}.
$$
\endproclaim
\par Hence we can define the linear map $ \beta: Fix(R)\rightarrow H^1(R) $ by the formula 
$$
 \beta(\mu) = F_\mu(R(z)) - R'(z)F_\mu(z).
$$  
In analogy with Kleinian group we  call $ \beta$ map as {\it Bers map }(see for example \cite{\bf Kra}).
\par Let $ A(S(R)) $ be the space of quadratic holomorphic integrable differentials on disconnected surface $ S(R).$ Let $ HD(S(R)) $ be the space of harmonic differentials on $ S(R)) $ that is locally in every chart every element $ \alpha \in HD(S(R)) $ has a form $ \alpha = \frac{\overline{\phi} d z^2}{\rho^2\vert dz\vert^2}, $ where $ \phi d z^2 \in A(S(R)) $ and $ \rho\vert d z\vert $ is Poincare metric. Let $ P :\{\overline{\Bbb C}\backslash{\big \{}\overline{\cup_iR^{-i}(Pc(R))}{\big \}}\rightarrow S(R) $ be the projection. Then the pullback $ P_*: HD(S(R))\rightarrow Fix(R) $ defines  a linear  injective map. \par The space $ HD(R) = P_*(A(S(R))) $ is called {\it the space of harmonic differentials.} For any element $ \alpha \in HD(R) $ the supporter $ supp(\alpha) \in F(R).$ Let $ HD(R) = P_*(HD(S(R)) \subset Fix(R) $ be subspace of harmonic differentials. Then $ dim(HD(R)) = dim(A(S(R))).$
\par Let $ J_R = Fix(R)_{\vert J(R)} $ be the space of invariant Beltrami differentials supported by Julia set.
\proclaim{Theorem 2} Let $ R $ be structurally stable rational map. Then $\beta: HD(R)\times J_R \rightarrow H^1(R) $ is an isomorphism.
\endproclaim
\par Note that in structurally unstable cases $ \beta $ restricted on $ HD(R)\times J_R $ is always injective.
\demo{Proof} $ R $ is structurally stable map hence $ dim(qc(R)) = dim(HD(R)\times J_R) = dim(H^1(R))$ $ = 2deg(R) - 2. $ If we can show that $ \beta $ is {\it onto} then we are done. 
\par Let us show that $\beta$ is "onto." 
\par Let $ G = z\frac{P_1Q - Q_1P}{Q^2} \in H^1(R) $ be any element. Let us consider the family of rational maps
$$
R_\lambda(z) = z\frac{P(z) + \lambda P_1(z)}{Q(z) + \lambda Q_1(z)}
$$
then $\frac{\partial R_\lambda(z)}{\partial\lambda}_{\vert\lambda = 0} = G(z) $ and for small $ \lambda $ maps $ R_\lambda \in K(R) $ (such that $ R $ is structurally stable). Let $ f_\lambda(z) $ be holomorphic family of qc-maps fixing points $  0, 1, \infty $ such that
$ R_\lambda = f_\lambda\circ R\circ f_\lambda^{-1}. $ The family of the complex dilatations $ \mu_\lambda(z) = \frac{\overline{\partial}f_\lambda(z)}{\partial f_\lambda(z)} \in Fix(R) $ forms meromorphic family of Beltrami differentials. If $ \mu_\lambda(z) = \lambda\mu_1(z) + \lambda^2\mu_2(z) + ..., $ where $ \mu_i(z) \in Fix(R). $  Then 
$$
\frac{\partial f_\lambda}{\partial\lambda}_{\vert\lambda = 0} = -\frac{z(z - 1)}{\pi}\iint\frac{\mu_1(\xi)}{\xi(\xi - 1)(\xi - z)} = F_{\mu_1}(z)
$$
and hence $ F_{\mu_1}(R(z)) - R'(z)F_{\mu_1}(z) = G(z). $
\par Now let $ \nu = {\mu_1}_{\mid F(R)}. $ Then we claim
\bigskip

{\bf Claim.} {\it There exists an element}  $ \alpha \in HD(S(R)) $
{\it such that} $ \beta(\alpha) = \beta(\nu). $
\bigskip
\par {\it Proof of claim.} We will use here qc-stuff  ( see for example books of I. Kra \cite{\bf Kra} S.L. Krushkal \cite{\bf Krush} papers C. McMullen and D. Sullivan \cite{\bf MM1} and \cite{\bf MS}). Let $ \omega $ be the Beltrami differential on $ S(R) $ generated by $ \nu $ (that is $ P_*(\omega) = \nu$). Let $ <\psi, \phi> $ be the Petersen scalar product on $ S(R), $ where $ \phi, \psi \in A(S(R)) $ and
$$
  <\psi, \phi>  = \iint_{S(R)} \rho^{-2}\overline{\psi}\phi,
$$
where $ \rho $ is hyperbolic metric on disconnected surface $ S(R). $
Then by (for example) Lemmas 8.1 and 8.2 of chapter III in \cite{\bf Kra} this scalar product defines a Hilbert space structure on $ A(S(R).
$ Then there exists an element $ \alpha^\prime \in HD(S(R)) $ such that equality 
$$ \iint_{S(R)}\omega\phi = \iint_{S(R)}\alpha^\prime\phi
$$
holds for all $ \phi \in A(S(R)).$
\par Now let $ A(O) $ be space of all holomorphic integrable over $ O $ functions, where $ O = \{\overline{\Bbb C}\backslash \overline{\cup_iR^{-i}(Pc(R))}\} \subset F(R). $ Then the pushforward operator $ P^*: A(O)\rightarrow A(S(R)) $ is dual to the pullback operator $ P_*.$ Hence element $ P_*(\alpha^\prime) $ satisfies the next condition
$$\iint_{O}\nu g = \iint_{O}P_*(\alpha^\prime) g
$$
for any $ g \in A(O). $
\par All above means that $ \iint P_*(\alpha)\gamma_a(z) = \iint\nu\gamma_a(z) $ for all $ \gamma_a(z) = \frac{a(a -1)}{z(z - 1)(z - a)}, a \in J(R). $ Hence the rational functions $ \beta(P_*(\alpha))(a) = \beta(\nu)(a) $ on $ J(R) $ and we have desired result with $ \alpha = P_*(\alpha^\prime). $ Claim and theorem are proved.
\enddemo
\subheading{Ruelle operator}
\proclaim{Definition} For a rational map $ R $ the Ruelle operator $ R^* $ is linear endomorphism of $ L_1({\Bbb C}) $ defined by the next formula
$$ R^*(\phi)(z) = \sum\phi(J_i(z))(J_i'(z))^2 = 
\sum_{y \in R^{-1}(z)}\frac{\phi(y)}{(R'(y))^2}.
$$
\endproclaim
Then we have the following  lemma.
\proclaim{Lemma 3} Let $ Y \subset \overline{\Bbb C} $ be completely invariant measurable subset respect to a rational map $ R.$ Then
\roster
\item $ R^* : L_1(Y) \rightarrow L_1{Y} $ is linear endomorphism "onto" with
$\Vert R^*\Vert_{L_1(Y)} \leq 1$;
\item The operator $ R_*(\phi) = \frac{\phi(R)(R^\prime)^2}{deg(R)}: L_1(Y)\rightarrow L_1(Y) $ is right inverse to $ R^*$;
\item  Beltrami operator $ B_R: L_\infty(Y)\rightarrow L_\infty(Y) $ is dual operator to $ R^*$;
\item If $ Y \subset \left \{\overline{\Bbb C}\backslash \overline{\cup_iR^{-i}(Pc(R))}\right \} $ is an open subset and let $ A(Y) \subset L_1(Y) $ be subset of holomorphic functions, then $ R^*(A(Y)) = A(Y). $
\endroster
\endproclaim
\demo{Proof} All items follow from definitions.
\enddemo
Let all critical points $ c_i $ be simple. Then there exists a decomposition $ \frac{1}{R'(z)} = \omega + \sum\frac{b_i}{z - c_i}, $ where $ \omega = \frac{1}{R'(\infty)} $ is multiplier of $ \infty $ and $ c_i $ are critical points and by residue theorem $ b_i = \frac{1}{R''(c_i)}.$ For $ i = 1, ..., 2deg(R) - 2 $ let $ h_i(z) = \frac{1}{R'(z)} - \frac{b_i}{z - c_i}. $ 
\proclaim{Lemma 4} For any rational map $ R $ with simple critical points we have
\roster
\item Let $ \gamma_a(z) = \frac{a(a - 1)}{z(z - 1)(z - a)} \in L_1({\overline{\Bbb C}}) $ where $ a \in {\Bbb C}\backslash \{0, 1\} $ is not a critical point, then 
$$ R^*(\gamma_a(z)) = \frac{\gamma_{R(a)}(z)}{R'(a)} + \sum_ib_i\gamma_a(c_i)\gamma_{R(c_i)}(z), $$
\par Let $ \tau_a(z) = \frac{1}{z - a} $ is a locally integrable function on $ \Bbb C, $ where $ a \in \Bbb C $ is not a critical point.
Then
$$
R^*(\tau_a(z)) = \frac{\tau_{R(a)}(z)}{R'(a)} + \sum_i b_i\tau_a(c_i)\tau_{R(c_i)}(z).
$$
\item if $ a = c_i $ is a critical point, then
$$ R^*(\gamma_a(z)) = (h_i(a) + b_i\frac{2c_i - 1}{c_i(c_i - 1)})\gamma_{R(a)}(z) + \sum_{j \neq i}b_j\gamma_a(c_j)\gamma_{R(c_j)}(z), $$
and 
$$
R^*(\tau_a(z)) = h_i(a)\tau_{R(a)}(z) + \sum_{j \neq i}b_j\tau_a(c_j)\tau_{R(c_j)}(z), 
$$
where coefficient $ h_i(a) + b_i\frac{2c_i - 1}{c_i(c_i - 1)} = \lim_{a\to c_i} \left(\frac{1}{R'(a)} + b_i\gamma_a(c_i)\right) $ 
\item If $ \mu \in Fix(R), $ then $ F_\mu(a) = -\frac{1}{\pi}\iint\mu\gamma_a(z) $ and 
$$ \beta(\mu)(a) = - R'(a)\sum b_iF_\mu(R(c_i))\gamma_a(c_i). $$
\endroster
\endproclaim
\demo{Proof} (1) Let $ \phi $ be any $ C^\infty $ function with compact support and $ supp(\phi) \in {\Bbb C}\backslash \{0, 1\}.$ Then we {\bf claim}
$$ \iint(\phi(R))_{\overline{z}}\gamma_a(z) = \iint\phi_{\overline{z}}\gamma_{R(a)}(z). $$
\par {\it Proof of the claim} By using decomposition $ \gamma_a(z) = \frac{a - 1}{z} - \frac{a}{z - 1} + \frac{1}{z - a} $ we obtain
$$\gather
\iint (\phi(R))_{\overline{z}}\gamma_a(z) = (a - 1)\iint\frac{(\phi(R))_{\overline{z}}}{z} - a\iint\frac{(\phi(R))_{\overline{z}}}{z - 1} + \iint\frac{(\phi(R))_{\overline{z}}}{z - a} = \\
= (a - 1)\phi(R(0)) - a\phi(R(1)) + \phi(R(a)) = *
\endgather
$$
such that $\phi(0) = \phi(1) = 0 $ and $ 0, 1 $ are fixed we have
$$ * = \phi(R(a)) = (R(a) - 1)\phi(0) - R(a)\phi(1) + \phi(R(a)) = \iint\phi_{\overline{z}}\gamma_{R(a)}(z). $$
Thus the claim is proved
\par Now let us show (1).
$$\gather
\iint\phi_{\overline{z}}R^*(\gamma_a(z)) = \iint\frac{\phi_{\overline{z}}(R){\overline{R'(z)}}}{R'(z)}\gamma_a(z) = \iint\frac{(\phi(R))_{\overline{z}}}{R'(z)}\gamma_a(z) =\\
= \omega\iint(\phi(R))_{\overline{z}}\gamma_a(z) + a(a - 1)\iint(\phi(R))_{\overline{z}}\sum \frac{b_i}{z(z - 1)(z - a)(z - c_i)} =
\endgather
$$
$$ \omega\iint(\phi(R))_{\overline{z}}\gamma_a(z) +  \sum\frac{a(a - 1)b_i}{a - c_i}\iint(\phi(R))_{\overline{z}}{\biggl (}\frac{1}{z(z -  1)(z - a)} - \frac{1}{z(z - 1)(z - c_i)}{\biggr )} = *
$$
then by the claim above we have
$$\gather * 
= \iint\phi_{\overline{z}}\gamma_{R(a)}(z){\biggl (}\omega + \sum\frac{b_i}{a - c_i}{\biggr )} - \sum\frac{b_i a(a - 1)}{c_i(c_i - 1)(a - c_i)}\iint\phi_{\overline{z}}\gamma_{R(c_i)}(z) = \\
= \iint\phi_{\overline{z}}\biggl(\frac{\gamma_{R(a)}}{R'(a)} + \sum b_i\gamma_a(c_i)\gamma_{R(c_i)}(z)\biggr). 
\endgather
$$
Let $ H(z) = R^*(\gamma_a(z)) - \frac{\gamma_{R(a)}}{R'(a)} - \sum b_i\gamma_a(c_i)\gamma_{R(c_i)}(z), $ then by above calculations we have $ \frac{\partial H}{\partial\overline{z}} = 0 $ (in sense of distributions) on $ {\Bbb C}\backslash\{0, 1\}.$  By the Weyl's lemma $ H $ is a holomorphic function on $ {\Bbb C}\backslash\{0, 1\}.$ and hence $ H(z) = C + \frac{A}{z} + \frac{B}{z - 1}. $ Besides the function $ H $ is integrable over $ \overline{\Bbb C}, $ hence $ H = 0.$ 
\par By argument above calculation of derivate of $ R^*(\tau_a(z)) $ (in sense of distributions) show that the function $ g(z ) = R^*(\tau_a(z)) - \frac{\tau_{R(a)}}{R'(a)} - \sum_i b_i\tau_a(c_i)\tau_{R(c_i)}(z) $ is holomorphic over $ \Bbb C $ and $ g(z)\to 0, $ for $ z\to \infty. $ Hence $ g(z) = 0. $
\par (2). Note that all integrals in calculations of the item (1) above depend on parameter $ a $ continuously. Hence we can consider limit when $ a \to c_i $ in formulas of the item (1). Then
$$
R^*(\gamma_{c_i}(z)) = \left(\lim_{a\to c_i}(\frac{1}{R'(a)} + b_i\gamma_a(c_i)\right)\gamma_{R(c_i)}(z) + \sum_{j \neq i}b_j\gamma_{c_i}{c_j}\gamma_{R(c_j)}(z)
$$
and 
$$
R^*(\tau_{c_i}(z)) = \left(\lim_{a\to c_i}(\frac{1}{R'(a)} + b_i\tau_a(c_i)\right)\tau_{R(c_i)}(z) + \sum_{j \neq i}b_j\tau_{c_i}{c_j}\tau_{R(c_j)}(z)
$$
We have $ \frac{1}{R'(a)} + b_i\gamma_a(c_i) = \omega + \sum_{j \neq i}\frac{b_j}{a - c_j} + \frac{b_i}{a - c_i} + \frac{b_i(a - 1)}{c_i} - \frac{b_i a}{c_i - 1} + \frac{b_i}{c_i - a}. $  Then in limit we are done.

\par(3) Let $ \mu \in Fix(R), $ then by item (1) we have
$$
\pi F_\mu(a) = -\iint\mu\gamma_a(z) = -\iint\mu R^*(\gamma_a(z)) = 
$$
$$ = -\frac{1}{R'(a)}\iint\mu\gamma_{R(a)}(z) - \sum_i b_i\gamma_a(c_i)\iint\mu\gamma_{R(c_i)}(z) =
$$
$$ = \pi\frac{F_\mu(R(a))}{R'(a)} + \pi\sum_i b_iF_\mu(R(c_i))\gamma_a(c_i).
$$
and hence
$$
\beta(\mu)(a) = F_\mu(a) - R'(a)F_\mu(a) = -R'(a)\sum_i b_iF_\mu(R(c_i))\gamma_a(c_i).\tag{*}
$$
So we are done. 
\enddemo

\proclaim{Remark 5} Lemma 4 above gives another coordinates for the spaces $ H^1(R) $ and $ HD(R)\times J_R. $ Namely the formula $ * $ above describes this isomorphism $ \beta^* : HD(R)\times J_R \rightarrow {\Bbb C}^{(2deg(R) - 2)} $ by
$$\beta^*(\mu) = (F_\mu(R(c_1)), ..., F_\mu(R(c_{2deg(R) - 2}))
$$
\endproclaim
\subheading{Formal relation of Ruelle-Poincare series}
\bigskip
\par Let again $ R \in X $ be a map with simple critical points.
\proclaim{Definition} Ruelle-Poincare series.
\roster
\item {\rm Backward  Ruelle-Poincare series.}
$$ RS(x, R, a) = \sum_{n = 0}^\infty(R^*)^n(\gamma_a)(x), $$ 
where $ a \in \overline{\Bbb C}\backslash \{0, 1\} $ is a parameter.
\item {\rm Forward Ruelle-Poincare series.}
$$ RP(x,R) = \sum_{n = 0}^\infty \frac{1}{(R^n)'(R(x))}.
$$
\item {\rm Modified Ruelle-Poincare series.} 
$$ A(x, R, a) = \sum_{n = 0}\frac{1}{(R^n)'(a)}\gamma_{R^n(a)}(x)
$$
\endroster
\endproclaim
Note that the Ruelle-Poincare series are a kind of generalizations of the Poincare series introduced by C. McMullen for rational maps (see \cite{\bf MM}). 
\par The next proposition gives a formal relation between Ruelle-Poincare series. 
\proclaim{Definition}  We denote the {\rm formal Cauchy product} of series $ A $ and $ B $ by $ A\otimes B $. Let us recall that if $ A = \sum_{i = 0} a_i $ and $ B = \sum_{i = 0} b_i, $ then $ C =  A\otimes B  = \sum_{i = 0} c_i, $ where $ c_i = \sum_{j = 0}^i a_jb_{i - j}.$ 
\endproclaim
\proclaim{Proposition 6} Let $ R \in X $ be rational map with simple critical points $ c_i. $ Let $ d_i = R(c_i) $ be critical values. Then we have the following formal relation between series
$$  
RS(z, R, a) = A(z, R, a) + \sum_i\frac{1}{R''(c_i)}A(c_i, R, a)\otimes RS(z, R, d_i).
$$
\endproclaim
Note that this proposition is part of Proposition A of preprint \cite{\bf Mak}.
\demo{Proof} By Lemma 4 we can calculate
$$ \gamma_a(z) = \gamma_a(z) $$
$$ R^*(\gamma_a(z)) = \frac{1}{R'(a)}\gamma_{R(a)}(z) + \sum_i b_i\gamma_a(c_i)\gamma_{d_i}(z)
$$
$$ (R^*)^2(\gamma_a(z)) = \frac{1}{(R^2)'(a)}\gamma_{R^2(a)}(z) + 
\sum_ib_i\biggl( \frac{\gamma_{R(a)}(c_i)}{R'(a)}\gamma_{d_i}(z) + 
\gamma_a(c_i)R^*(\gamma_{d_i}(z))\biggr)
$$
$$ . . . $$
$$\multline 
(R^*)^n(\gamma_a(z)) = \frac{1}{(R^n)^\prime(a)}\gamma_{R^n(a)}(z) +\\
+ \sum_i b_i{\biggl (}\frac{\gamma_{R^{n - 1}(a)}(c_i)}
{(R^{n - 1})^\prime(a)}\gamma_{d_i}(z) + 
\frac{\gamma_{R^{n - 2}(a)}(c_i)}{(R^{n - 2})^\prime(a)}R^*(\gamma_{d_i}(z)) 
+ ... + 
\gamma_a(c_i)(R^*)^{n - 1}(\gamma_{d_i}(z))\biggr).
\endmultline
$$
Summation of columns gives the desired expression.
\enddemo
For a map $ R $ and  $ \vert x\vert < 1 $ define the following formal series
$$ RS(x, z, R, a) = \sum_{i = 0} x^i(R^*)^i(\gamma_a(z)) \text{ and }
$$
$$ A(x, z, R, a) = \sum_{i = 0}\frac{x^i}{(R^i)^\prime(a)}\gamma_{R^i(a)}(z).
$$
Then we have the following lemma.
\proclaim{Lemma 7} Let $ R $ be rational map and $ 0, 1, \infty $ be 
fixed points. Then
\roster 
\item $ RS(x, z, R, a) \in L_1(\overline{\Bbb C}) $ for any 
$ a \in {\Bbb C} $ and $ \vert x\vert < 1.$
\item Assume that a point $ a \in {\Bbb C} $ is summable for $ R. $ Then  the 
function $ A(x, z, R, a)$ $\in$ $ L_1(\overline{\Bbb C}) $ for all $ \vert x\vert \leq 1. $ Moreover
$$ \lim_{x\to 1}\Vert A(x, z, R, a) - A(1, z, R, a)\Vert = 0.$$
\endroster
\endproclaim
\demo{Proof} (1) The norm $ \Vert R^*\Vert \leq 1 $ hence
$$ \iint\vert RS(x, z, R, a)\vert \leq \Vert\gamma_a(z)\Vert\sum_{i = 0}^{\infty}\vert x\vert^i = \frac{\Vert\gamma_a(z)\Vert}{1 - \vert x\vert}.
$$
\par (2) Due to classical theorems (see for example books of Kra \cite{\bf Kra} and Vekua \cite{\bf Vek}) there exists a constant $ M < \infty $ which does not depend on $ a $ such that we have the following estimate $ \Vert\gamma_a(z)\Vert \leq M\vert a\vert\vert\ln\vert a\vert\vert.$ 
\par Then we have
$$ \iint\vert A(x, z, R, a)\vert \leq M\sum{\biggl |}\frac{x^i\vert R^i(a)\vert\vert\ln\vert R^i(a)\vert\vert}{(R^i)'(a)}{\biggr |} < \infty.
$$
\par Now hence we have the following estimate
$$ \lim_{x\to 1}\Vert A(x, z, R, a) - A(1, z, R, a)\Vert \leq M\sum{\biggl |}\frac{(x^i - 1)\vert R^i(a)\vert\vert\ln\vert R^i(a)\vert\vert}{(R^i)'(a)}{\biggr |}
$$
Let us show that $ \lim_{x\to 1}\sum{\biggl |}\frac{(x^i - 1)\vert R^i(a)\vert\vert\ln\vert R^i(a)\vert\vert}{(R^i)'(a)}{\biggr |} = 0. $ 
\par To do it let $ \epsilon > 0 $ be any fixed number. Choose $ N $ such that $ 2\sum_{i \geq N}{\biggl |}\frac{\vert R^i(a)\vert\vert\ln\vert R^i(a)\vert\vert}{(R^i)'(a)}{\biggr |} \leq \frac{\epsilon}{2}. $ Let $ \delta $ be such number that for all $x, \vert 1 - x\vert <\delta $  we have $ \frac{\vert 1 - x^N\vert}{C} < \frac{\epsilon}{2}, $ where $ C = \sum{\biggl |}\frac{\vert R^i(a)\vert\vert\ln\vert R^i(a)\vert\vert}{(R^i)'(a)}{\biggr |}, $ then we have the following estimate
$$\multline  
\sum{\biggl |}\frac{(x^i - 1)\vert R^i(a)\vert\vert\ln\vert R^i(a)\vert\vert}{(R^i)'(a)}{\biggr |} \leq \sum_{i < N}{\biggl |}\frac{(x^i - 1)\vert R^i(a)\vert\vert\ln\vert R^i(a)\vert\vert}{(R^i)'(a)}{\biggr |} +\\
 \sum_{i \geq N}{\biggl |}\frac{(x^i - 1)\vert R^i(a)\vert\vert\ln\vert R^i(a)\vert\vert}{(R^i)'(a)}{\biggr |} \leq \frac{\epsilon}{2} + 2\sum_{i \geq N}{\biggl |}\frac{\vert R^i(a)\vert\vert\ln\vert R^i(a)\vert\vert}{(R^i)'(a)}{\biggr |} \leq \epsilon.
\endmultline
$$
Thus the lemma is finished.
\enddemo
\proclaim{Lemma 8} If $ a \in {\Bbb C} $ is a summable point for a structurally stable rational map $ R. $ Then for any fixed $ \vert x\vert \leq 1 $ the function $ A(x, z, R, a) $ is meromorphic function (respect to variable $ z $) which is finite in all  critical points. Moreover
$$ \lim_{x\to 1} A(x, c, R, a) = A(1, c, R, a).
$$
\endproclaim
\demo{Proof}  Let $ c  $ be a critical point. If $ c \notin \overline{\cup_nR^n(a)} $ and $ d $ is the distance between $ c $ and $ \overline{\cup_nR^n(a)}, $  then
$$
{\biggl |}A(x, c, R, a){\biggr |} \leq \frac{M}{\vert c(c - 1)d\vert}\sum_i{\biggl |}\frac{x^i}{(R^i)'(a)}{\biggr |} \leq \frac{M}{\vert c(c - 1)d\vert}\sum\frac{1}{\vert (R^i)'(a)\vert} < \infty
$$
and arguments of lemma 7 above complete the proof.
\par Now assume that $ c \in \overline{\cup_nR^n(a).} $  Let $ U_\epsilon $ be $\epsilon-$neighborhood of $ c. $ Let points $ R^{n_i}(a) \in U_\epsilon. $ Then by arguments above and $ c \neq R^n(a) $ for all $ n \geq 0 $ we have to estimate the following expression
$$        \sum_i\frac{x^{n_i}}{(R^{n_i})'(a)}\gamma_{R^{n_i}(a)}(c).
$$
By using equality $ R'(z) = (z - c)R''(c) + O(\vert z- c\vert^2 $ for $ z \in U_\epsilon $ we obtain
$$ \frac{1}{\vert R^{n_i} - c\vert} \leq \frac{\vert R''(c)\vert + O(\vert R^{n_i} - c\vert)}{\vert R'(R^{n_i}(a))\vert} \leq C\frac{1}{\vert R'(R^{n_i}(a))\vert}
$$
hence
$$
{\biggl |}\gamma_{R^{n_i}(a)}(c){\biggr |} \leq C\frac{1}{\vert R'(R^{n_i}(a))\vert}\frac{\vert R^{n_i}(a)\vert\vert R^{n_i}(a) - 1\vert}{\vert c(c - 1)} \leq C_1 \frac{1}{\vert R'(R^{n_i}(a))\vert}
$$
where $ C $ and $ C_1 $ are constant depending only on $\epsilon $ and the point $ c. $
As result for all $ \vert x\vert \leq 1 $ we have
$$ {\biggl |}\sum_i\frac{x^{n_i}}{(R^{n_i})'(a)}\gamma_{R^{n_i}(a)}(c){\biggr |} \leq C_1\sum_i\frac{\vert x\vert^{n_i}}{\vert (R^{n_i + 1})'\vert} < \infty.
$$
	\par The arguments of Lemma 7 above show that
$$ \lim_{x\to 1} A(x, c, R, a) = A(1, c, R, a).
$$
The lemma is finished.
\enddemo    
\proclaim{Corollary 9} Let $ R $ be a rational map with simple critical points. Assume a point $ a \in {\Bbb C}\backslash \{0, 1\} $ is summable for $ R. $  Then  respect to variable $ x $ the following function equality holds.
$$ RS(x, z, R, a) = A(x, z, R, a) + x\sum_i b_i RS(x, z, R, d_i)\cdot A(x, c_i, R, a).\tag{*}
$$
\endproclaim

\demo{Proof} Lemmas 2-4 and properties of Cauchy product give the desired equality.
\enddemo
 \par Now we are ready to prove Theorem A. 
\proclaim{Theorem 10} Let $ R \in X $ be a rational map, then $ dim(HD(R)\times J_R) < 2d - 2.$
\endproclaim
\demo{Proof} Let $ c_1 \in J(R) $ be a summable critical point.  
\par Assume $ dim(HD(R)\times J_R) = 2d - 2,$ then $ R $ is structurally stable map. Let us consider relation between Ruelle-Poincare series due by corollary 9 for $ a = R(c_1) = d_1. $ Let $\mu \in HD(R)\times J_R $ be any element. Let us integrate $ \mu $ with relation $ * $ like follows
$$ \iint\mu RS(x, z, R, d_1) =\iint\mu A(x, z, R, d_1) + x\sum_i b_i A(x, c_i, R, d_1)\cdot\iint\mu RS(x, z, R, d_i).
$$
Then using the invariance of $ \mu $ we obtain
$$
F_\mu(d_1)\sum_i x^i = \sum_i\frac{x^i F_\mu(R^i(d_1))}{(R^i)'(d_1)} + x\sum_j b_j F_\mu(d_j){\biggl (}\sum_ix^i{\biggr )} A(x, c_j, R, d_i).
$$
and hence
$$\frac{F_\mu(d_1)}{(1 - x)} = \sum_i\frac{x^i F_\mu(R^i(d_1))}{(R^i)'(d_1)} + \frac{x}{(1 - x)}\sum_j b_j F_\mu(d_j)A(x, c_j, R, d_i).\tag{2}
$$
Let $ C_j = \lim_{x\to 1}A(x, c_i, R, d_1) $ by Lemmas 3 and 4 $ C_j $ are good defined constants which do not depend on $ \mu.$ By using the following estimate
$$ {\biggl |}(1-x)\sum_ix^i\frac{F_\mu(R^i(d_1))}{(R^i)'(d_1)}{\biggr |} \leq M_\mu\vert 1 - x\vert\sum_i\frac{1}{\vert (R^i)'(d_1)\vert}\leq C_\mu\vert 1 - x\vert.
$$
and taking limit for $ x\to 1 $ in $(2) $ we obtain the equation.
$$ F_\mu(d_1) = \sum_i b_iC_iF_\mu(d_i) \text{ or }  
F_\mu(d_1)(1 - b_1C_1) = \sum_{i \geq 2} b_iC_iF_\mu(d_i).\tag{3}
$$
Here the coefficients $ C_i $ do not depend on $ \mu.$
	\par Now we  call relation (3) as {\it trivial relation } iff
$$ C_1 = \frac{1}{b_1} \text{ and } C_i = 0  \text{ for }  i > 1.
$$
Here we need the following proposition
\enddemo
\proclaim{Proposition 11} If relation (3) above is non-trivial, then $ R $ is structurally unstable.
\endproclaim
\demo{Proof} By Remark 5 the operator $ \beta $ induces an isomorphism
$$ \beta^*: HD(R)\times J_R\rightarrow {\Bbb C}^{(2d - 2)}.
$$
If relation (3) above is non-trivial, then this relation gives non-trivial equation on image of $ \beta^* $ and $ \Im(\beta^*) $ is subset of the set of solutions of this equation. Hence $ dim(HD(R)\times J_R) = dim(\Im(\beta^*)) < 2d - 2 $ and we have desired conclusion. 
\enddemo
	\proclaim{Proposition 12}
\roster
\item If relation (3) is trivial, then
$$
\sum_{i = 0}\frac{1}{(R^i)'(d_1)} = \sum_{i = 0} \frac{R^i(d_1)}{(R^i)'(d_1)} = 0
$$
\item If $ \phi(z) = A(z, R, d_1)\neq 0 $ identically on $ Y = {\Bbb C} \backslash X_{c_1}, $ then relation (3) is non - trivial. 
\endroster
\endproclaim
\par Here we use the idea of G. Levin \cite{\bf Lev} of consideration of function $ \phi(z) = A(z, R, d_1) $ itself. Levin observes that (in case deg(R) = 2) that the function $ A(z, R, d_1) $ is a fixed point for $ R^*. $
\demo{Proof} Assume that relation (3) is trivial relation, that is
$$
C_1 = A(c_1, R, d_1) = \phi(c_1) = \frac{1}{b_1} \text{ and }  \,C_i = A(c_i, R, d_1) = \phi(c_i) = 0 \text{ for } i > 1.
$$
Then we {\bf claim.} {\it Under conditions above we have}
$$ R^*(\phi(z)) = \phi(z).
$$
\demo{Proof of the claim} By lemma 4 calculations show
$$ R^*(\phi(z)) = \phi(z) - \gamma_{d_1}(z) + \sum_i^{2d - 2} b_i\gamma_{d_i}A(c_i, R, d_1)
$$
but under assumptions $ \phi(c_i) = C_i = 0 $ for $ i > 1 $ and $ \phi(c_1) = \frac{1}{b_1} $ and hence $ R^*(\phi(z) = \phi(z).$
The claim is finished.
\enddemo
\par Now let us check (1). Let $ A = \left (\sum_{n = 0}\frac{R^n(d_1) }{(R^n)'(d_1)}\right ) $ and $ B = \left (\sum_{n = 0}\frac{ 1}{(R^n)'(d_1)}\right ). $ Then we have
$$
\phi(z) = \frac{1}{z}\left (\sum_{n = 0}\frac{R^n(d_1) - 1}{(R^n)'(d_1)}\right ) - \frac{1}{z - 1}\left (\sum_{n = 0}\frac{1}{(R^n)'(d_1)}\right ) + \psi(z) = \frac{A - B}{z} - \frac{A}{z - 1} + \psi(z),
$$
where $ \psi(z) = \sum_i\frac{\tau_{R^i(d_1)}(z)}{(R^i)'(d_1)}. $
Hence
$$ 
R^*(\phi(z)) = (A - B)R^*(\tau_0(z)) - B R^*(\tau_1(z)) + \sum_i \frac{R^*(\tau_{R(d_1)}(z))}{(R^i)'(d_1)} = \ast
$$
by the lemma 4 (calculation for $\tau_a(z)$) we have
$$ 
\ast 
 = (A - B)\left(\frac{\tau_0(z)}{R'(0)} + \sum_j b_j\tau_0(c_j)\tau_{R(c_j)}(z)\right ) - B\left (\frac{\tau_1(z)}{R'(1)} + \sum_j b_j\tau_1(c_j)\tau_{R(c_j)}(z)\right) + \psi(z) -
$$
$$\multline
 - \tau_{d_1}(z) +\sum_j b_j\tau_{R(c_j)}(z)\psi(c_j) = 
\sum_j b_j\tau_{R(c_j)}(z)\biggl((A - B)\tau_0(c_j) - B\tau_1(c_j) + \psi(c_j)\biggr) +\\
+ \biggl(\psi(z) + (A - B)\tau_0(z) - B\tau_1(z)\biggr) + 
\biggl((A - B)\tau_0(z)\left(\frac{1}{R'(0)} - 1\biggr) -  B\tau_1(z)\left(\frac{1}{R'(1)} - 1\right)\right ) -\\
- \tau_{d_1}(z) = \phi(z) - \tau_{d_1}(z) + \sum_j b_j\tau_{R(c_j)}(z)\phi(c_j)  +\\
+ \left ((A - B)\tau_0(z)\left(\frac{1}{R'(0)} - 1\right) - B\tau_1(z)\left(\frac{1}{R'(1)} - 1\right)\right ) = \ast\ast
\endmultline
$$
Such as relation (3) is trivial we have
$$
\ast\ast = \phi(z) + (A - B)\tau_0(z)\left(\frac{1}{R'(0)} - 1\right) - B\tau_1(z)\left(\frac{1}{R'(1)} - 1\right) = \text{ by claim above } = \phi(z).
$$
Hence $ A = B = 0. $ We are finished (1). \hfill\newline
2) We need the following lemma

\proclaim{Lemma 13} Under conditions above 
the measure $ \sigma(A) = \iint_A\vert\phi(z)\vert $
is non-negative invariant absolutely continue probability measure, where $ A\subset \overline{\Bbb C} $ is measurable subset.
\endproclaim 
\demo{Proof} In this proof we use one observation of A. Epstein (see \cite{\bf Eps}.)  Invariance the measure $ \sigma $ means that 
$$ \vert\phi\vert = \vert R^*\vert\vert\phi\vert = \sum_i^{deg(R)}\vert \phi(J_i)\vert\vert J_i'\vert^2
$$
where operator $\vert R^*\vert $ is called {\it the modulus of Ruelle operator}.
\par Let us show invariance.  For any fixed $ 1 \leq j\leq deg(R) $ claim above and calculations below show. 
$$ \Vert\phi(z)\Vert = \Vert R^*(\phi)\Vert = \iint\vert\sum_i\phi(J_i)(J_i')^2\vert \leq \iint\vert\phi(J_j)(J_j')^2\vert + \iint\vert\sum_{i\neq j}\phi(J_i)(J_i')^2\vert \leq\Vert\phi\Vert.
$$
In expression above we have the right part is equal to left last part. Hence all signs inequality are indeed equality. Let us define $ \alpha_j(z) = \phi(J_j(z))(J_j(z)')^2 $ and $ \beta_j(z) = \phi(z) - \alpha_j(z) = \sum_{i\neq j}\phi(J_i(z))(J_i(z)')^2. $ Then by above we have
$$ \iint\vert\alpha_j + \beta_j\vert = \iint\vert\alpha_j\vert + \iint\vert\beta_j\vert.
$$
We deduce that $ \vert\alpha_j + \beta_j\vert = \vert\alpha_j\vert + \vert\beta_j\vert $ almost everywhere. Indeed otherwise let $ A = \{z, \vert\alpha_j(z) + \beta_j(z)\vert < \vert\alpha_j(z)\vert + \vert\beta_j(z)\vert\}$ with $ m(A) > 0, $ where $ m $ is the Lebesgue measure. Then $ \Vert\phi\Vert = \iint_A\vert\phi\vert + \iint_{\overline{\Bbb C}\backslash A}\vert\phi\vert <\Vert\phi\Vert. $Contradiction. By induction we finish lemma 13.
\enddemo
\proclaim{Corollary 14} If $ \phi(z) \neq 0 $ identically on $ Y_a, $
then $ J(R) = \overline{\Bbb C} $ and $ \frac{\overline{\phi}}{\vert\phi\vert} $ defines an invariant Beltrami differential.
\endproclaim
\demo{Proof} Recall that a measurable positive Lebesgue measure set $ A \subset \overline{\Bbb C} $ is {\it wandering} for  a map $ R $ iff $ m(R^{-n}(A)\cap R^{-k}(A)) = 0 $ for $ k \neq n. $ If $ \lambda $ is non-negative invariant probability measure, then $ \lambda(A) = 0 $ for any wandering set $ A.$ Then we have
\roster
\item If $ F(R) \neq\emptyset, $ then by lemma 13 and argument above $ \sigma(F(R)) = 0 $ and hence $ \phi = 0 $ identically on $ Y_a.$
\item If $ J(R)  = \overline{\Bbb C}, $ then either $ \phi = 0 $ identically on $ Y_a $ or $ \phi \neq 0 $ identically on every component of $ Y_a .$
\endroster
	In the last case let us show that
$ \mu = \frac{\overline{\phi}}{\vert\phi\vert} \in Fix. $ Indeed  in notations of lemma 12 let us show that  
$$\frac{\beta_j}{\alpha_j} = k_j \geq 0 \text{ is non negative function. }
$$
Really we have $ \vert 1 + \frac{\beta_j}{\alpha_j}\vert = 1 + {\biggl |}\frac{\beta_j}{\alpha_j}{\biggr |}, $ then if $ \frac{\beta_j}{\alpha_j} = \gamma^j_1 + i\gamma^j_2 $ we have
$$ {\biggl (}1 + (\gamma^j_1){\biggr )}^2 + (\gamma^j_2)^2 = {\biggl (}1 + \sqrt{(\gamma_1^j)^2 +(\gamma_2^j)^2}{\biggr )}^2 = 1 + (\gamma_1^j)^2 +(\gamma_2^j)^2 + 2\sqrt{(\gamma_1^j)^2 +(\gamma_2^j)^2}.
$$
Hence $ \gamma_2^j = 0  $ and $ \frac{\alpha_j}{\beta_j} = \gamma_1^j $ is a real-valued function but $\frac{\alpha_j}{\beta_j}$ is meromorphic function. So $ \gamma_1^j = k_j $ is constant on every connected component of $ Y_a $ and condition $ \vert 1 + k_j\vert = 1 +\vert k_j\vert $ shows $ k_j \geq 0.$
\par Then for any $ 1 \leq j \leq deg(R) $ we have
$$ \frac{R^*(\phi)}{\alpha_j} = \frac{\alpha_j +\beta_j}{\alpha_j} = 1 + k_j \geq 1
$$
and
$$ \phi = R^*(\phi) = (1 + k_j)\phi(J_j)(J_j')^2. $$
That means for any $ j $
$$ \mu(z) = \frac{(1 + k_j)\overline{\phi(J_j)(J_j')^2}}{\vert(1 + k_j)\phi(J_j)(J_j')^2\vert}= \mu(J_j)\frac{\overline{J_j'}}{J_j'}
$$
and we obtain
$$ \mu(R)\frac{\overline{R'}}{R'} = \mu.
$$
Corollary is complete.
\enddemo
\par To finish Proposition 12 we use the McMullen result (lemma 3.16, \cite{\bf MM1}) which states that if  $ \frac{\overline{\phi}}{\vert\phi\vert} \in Fix $ for integrable holomorphic function $ \phi(z) \neq 0 $ on a domain $ U \subset J(R). $ Then  $ R $ is {\it double covered by an integral torus endomorphism} and particularly it is unstable. Proposition 12 is finished. 
\proclaim{Corollary A}  Let $ R $ be a rational map with summable critical point $ c \in J(R).$ If the sum $ \sum_{i = 0}\frac{1}{(R^i)\prime(R(c))} \neq 0, $ then $ R $ is unstable map.
\endproclaim
\demo{Proof} By conditions and proposition 12 we have the relation (3) can not be trivial. Application proposition 11 finished proof of this corollary.
\enddemo
\par To finish theorem 10 we need the next proposition.
 \proclaim{Proposition 15} Let $ a_i \in {\Bbb C}, a_i\neq a_j, $ for $ i\neq j $ be points such that $ Z = \overline{\cup_i a_i} \subset {\Bbb C} $ is a compact set. Let  $ b_i \neq 0 $ be complex numbers such that the series $ \sum b_i $ is absolutely convergent. Then the function $ l(z) = \sum_i\frac{b_i}{z - a_i} \neq 0 $ identically on $ Y = {\Bbb C} \backslash Z $ in the following cases
\roster
\item the set $ Z $ has zero Lebesgue measure
\item if diameters of components of $ {\Bbb C} \backslash Z $ uniformly bounded below from zero and
\item If $ O_j $ denote the components of $ Y, $ then $ \cup_i a_i \in \cup_j{\partial O_j}.$
\endroster
\endproclaim
\demo{Proof}  Assume that $ l(z) = 0 $ on $ Y.$ Let us calculate derivative $ \overline{\partial}l $ in sense of distribution, then $ \omega = \overline{\partial}l = \sum_i b_i\delta_{a_i} $ and by standard arguments 
$$ l(z) = - \int\frac{d\omega(\xi)}{\xi - z}.
$$
Such as $ a_i \neq a_j $ for $ i \neq j, $  then measure $ \omega =0 $ iff all coefficients $ b_i = 0. $ 

Let us check (1). Otherwise in this case we have that the function $ l $ is locally integrable and $ l = 0 $ almost everywhere and hence $ \omega = \overline{\partial}l = 0 $ in sense of distributions and hence $ \omega = 0 $ as a functional on space of all continuous functions on $ Z $ that is contradictions with arguments above.

2)   Assume that $ l = 0 $ identically out of $ Z.$ Let $ R(Z) \subset C(Z) $ denote the algebra of all uniform limits of rational functions with poles out of $ Z $ in $\sup-$topology, here $ C(Z) $ as usually denotes the space of all continuous functions on $ Z $ with $ \sup-$norm. Then measure $ \omega $ denote a lineal functional on $ R(Z).$ 
The items (2) and (3) are based on the generalized Mergelyan theorem (see \cite{\bf Gam} thm. 10.4) which states {\it If diameters of all components of $ {\Bbb C} \backslash Z $ are bounded uniformly below from 0, then  every continuous function holomorphic on interior of $ Z $ belongs to $ R(Z). $}

Let us show that $ \omega $ annihilates the space $ R(Z). $ Indeed let $ r(z) \in R(Z) $ be rational function and $ \gamma $ enclosing $ Z $ close enough to $ Z $ such that $ r(z) $ does not have poles in interior of $ \gamma.$ Then such that $ l = 0 $ out of $ Z $ we only apply Fubini's theorem
$$ 
\int r(z)d\omega(z) = \int d\omega(z)\frac{1}{2\pi i}\int_{\gamma}\frac{r(\xi)d\xi}{\xi - z} = \frac{1}{2\pi i}\int_{\gamma}r(\xi)d\xi\int\frac{d\omega(z)}{\xi - z} = \frac{1}{2\pi i}\int_{\gamma}r(\xi)l(\xi)d\xi = 0.
$$
Then by generalized Mergelyan theorem we have $ R(Z) = C(Z) $ and $ \omega = 0. $ Contradiction.

\par Now let us check (3). We {\bf claim} {\it that $ l = 0 $ almost everywhere on} $ \cup_i\partial O_i. $

\par {\it Proof of the claim.}  Let $ E \subset \cup_i\partial O_i $ be any measurable subset with positive Lebesgue measure. Then the function $ F_E(z) = \iint_E\frac{dm(\xi)}{\xi - z} $ is continuous on $ {\Bbb C} \backslash \cup_i O_i $ and is holomorphic onto interior of $ {\Bbb C} \backslash  \cup_i O_i. $ Again by generalized Mergelyan theorem  $ F_E(z) $ can be approximated on $ {\Bbb C} \backslash \cup_i O_i $ by functions from $ R({\Bbb C} \backslash  \cup_i O_i) $ and hence by arguments above and by assumption we have
$ \int F_E(z)d\omega(z) = 0. $ But again application of Fubini's theorem gives
$$ 0 = \int F_E(z)d\omega(z) = \int d\omega(z)\iint_E\frac{dm(\xi)}{\xi - z} = \iint_E dm(\xi)\int d\omega(z)\frac{1}{\xi - z} = \iint l(\xi)d(m(\xi).
$$
Hence for any measurable $ E \subset \cup_i\partial O_i $ we have $ \iint_E l(z) = 0. $ The claim is proved.
Now for any component $ O \in Y $ and  any measurable $ E \subset \partial O $ we have $ \iint_E l(z) = 0. $ By assumption $ l = 0 $ almost everywhere on $ \Bbb C.$ Contradiction. Proposition is proved.
\enddemo

\proclaim{Proposition 16} Under condition of theorem 10 the function $ A(z, R, d_1) = \phi(z) \neq 0 $ identically on $ Y $ in the following cases 
\roster
\item if $ c_1 \notin X_{c_1}, $
\item if diameters of components of $ Y $ are uniformly bounded below from 0, 
\item  If $ m\left (X_{c_1}\right ) = 0, $ where $ m $ is the Lebesgue measure on $ {\Bbb C}, $
\item if $ X_{c_1} \subset \cup_i \partial D_i, $ where $ D_i $ are components of Fatou set. 
\endroster
\endproclaim
\proclaim{Remark} A. Eremenko (personal communication) can prove this proposition for polynomials. He uses harmonic functions stuff. But we believe that  our arguments (proposition 15) are more simple.
\endproclaim
\demo{Proof} 
1)If $ R $ is structurally stable then relation (3) is trivial. Then particularly $ \phi(c_1) = \frac{1}{b_1} \neq 0 $ and we are done case 1).\hfill\newline
\par  First assume that the set $ X_{c_1} $ is bounded. Then\hfill\newline
By proposition 12 we have that $ \phi(z) = \psi(z). $ Hence by application proposition 15 we are finished last 3 cases. 
\par Now let $ X_{c_1} $ be unbounded. Let $ y \in {\Bbb C} $ be a point such that the point $ 1 - y \in Y, $ then the map $ g(z) = \frac{y z}{z + y - 1} $ maps $ X_{c_1} $ into $ \Bbb C. $ Let us consider the function $ f(z) = \sum\frac{1}{(R^i)'(d_1)(z - g(R^i(d_1))}, $ then by proposition 15 $ f(z) \neq 0 $ identically on $ g(Y). $
\par Now we {\bf Claim} that {\it Under condition of theorem 10 $ f(g(z)) = (z + y - 1)^2\phi(z)$.}
\par {\it Proof of claim.} We have for any $ n $
$$
\frac{1}{g(z) - g(R^n(d_1))} = \frac{(z + y - 1)(R^n(d_1) + y - 1)}{y(y - 1)(z - R^n(d_1))} = \frac{1}{y(y - 1)}\left(\frac{(z + y - 1)^2}{z - R^n(d_1)} + 1 - y - z\right),
$$
then 
$$
\multline
f(g(z)) = \sum \frac{1}{(R^i)'(d_1)(g(z) - g(R^i(d_1))} =\\ 
= \frac{1}{y(y - 1)}\biggl((1 - y - z)\sum \frac{1}{(R^i)'(d_1)} + (z + y - 1)^2\sum\frac{1}{(R^i)'(d_1)(z - g(R^i(d_1))}\biggr) = \ast
\endmultline
$$
and by proposition 12 we have
$$
\ast = (z + y - 1)^2\phi(z)
$$
Claim is finished.
\par If $ \phi = 0 $ on $ Y, $ then $ f = 0 $ on $ g(Y). $ Contradiction with proposition 15. Proposition 16 is proved.
\enddemo

\par To finish theorem 10 we use propositions 11 and 16.
\enddemo

\heading{\bf Proof of Theorem B}
\endheading
\par Let us formulate one more the {\bf main observation} of this paper.
\proclaim{Main observation} Each summable critical point $ c_1 $ with bounded forward orbit gives an equation on image of the operator $ \beta. $ Namely for all $ \mu \in HD(R)\times J_R $ we have
$$ F_\mu(R(c_1))(1 - b_1 C_1) = \sum_{i > 1} b_iC_iF_\mu(R(c_i))
$$
and hence
$$ \beta(\mu)(a) = -R'(a)\sum_{i \geq 2}b_iF_\mu(R(c_i)){\biggl (}\gamma_a(c_i) +\frac{C_ib_1}{1 - b_1 C_1}\gamma_a(c_1){\biggr )}
$$
and hence $ \dim((\beta(HD(R)\times J_R)) \leq 2deg(R) - 3. $ 
\endproclaim
	By arguments above we obtain that image  $ \beta(HD(R)\times J_R) $ belongs to common solution of the following system of equations
$$  F_\mu(R(c_i))(1 - b_i A(c_i, R, R(c_i)) = \sum_{j \neq i} b_jF_\mu(R(c_j))A(c_j, R, R(c_i)),
$$
for all $ c_i \in J(R), i =1, ..., k. $ Hence if this system is linearly independent, then 
$$ dim(\beta(HD(R)\times J_R)) = dim(qc(R)) - \#\{\text{summable critical points on }  J(R)\} 
$$
 and by using injectivity of $ \beta $ we have $ J_R = \emptyset. $ Theorem is complete.
\par Now assume that  system above is linearly dependent. That means that there exist numbers $ B_i $ such that the function
$$ \alpha(z) = \sum B_i A(z, R, R(c_i)) 
$$
is fixed points for Ruelle operator $ R^*. $ By arguments of theorem A we obtain that $ \alpha(z) = 0 $ almost everywhere on $ {\Bbb C}. $ That means that measure 
$$ \frac{\partial\alpha}{\partial\overline{z}} = \sum_iB_i\sum_n\frac{\delta_{R^n(R(c_i))}}{(R^n)'(R(c_i))} = 0.
$$
Then we obtain contradiction with condition (2) in definition of $ W. $
Now by using proposition 13 we are finish theorem B.

\enddocument
\end